\documentclass[
	a4paper, 
	8pt, 
]{LTJournalArticle}

\usepackage{etoolbox}
\providetoggle{preprint}
\settoggle{preprint}{true}
\providetoggle{singlecolalg}
\settoggle{singlecolalg}{false}

\usepackage{authblk}
\usepackage{graphicx}
\usepackage{amsthm}

\usepackage[utf8]{inputenc}
\usepackage[textsize=scriptsize,shadow,disable,loadshadowlibrary]{todonotes}
\setuptodonotes{inline}
\usepackage{url}
\usepackage{bm}
\usepackage{enumitem}
\usepackage{verbatim}

\usepackage{hyperref}
\hypersetup{
	colorlinks = true,
	linkcolor=black,
	citecolor={blue!50!black},
	urlcolor={blue!80!black},
	anchorcolor = blue
}
\hypersetup{linkcolor=black}

\usepackage{amsmath}
\usepackage{amssymb}
\usepackage{mathtools}
\usepackage{siunitx}[=v2]
\usepackage{booktabs}
\usepackage[super]{nth} 

\usepackage[linesnumbered,ruled,vlined]{algorithm2e}
\usepackage{setspace}
\definecolor{commentgr}{rgb}{0,0.7,0}

\SetCommentSty{commfont}
\SetKwComment{Comment}{\#}{}
\RestyleAlgo{ruled}
\SetKwRepeat{Do}{do}{while}
\SetKw{KwBreak}{break}
\SetKwBlock{With}{With}{end}
\SetKwBlock{Forpar}{for parallel}{end}

\DeclareMathOperator*{\argmin}{arg\,min}

\newcommand{\norm}[1]{\left\lVert#1\right\rVert}
\newcommand{\abs}[1]{\left\lvert#1\right\rvert}

\newcommand{\proj}[0]{\text{proj}} 
\newcommand{\diff}[1]{\ensuremath{\operatorname{d}\!{#1}}}

\newcommand{\vx}[0]{\ensuremath{\boldsymbol{x}}}

\newcommand{\vv}[0]{\ensuremath{\boldsymbol{v}}}

\newcommand{\vy}[0]{\ensuremath{\boldsymbol{y}}}
\newcommand{\balpha}[0]{\ensuremath{\boldsymbol{\alpha}}}

\newcommand{\bPhi}[0]{\ensuremath{\boldsymbol{\Phi}}}

\newcommand{\tD}{\mathbf{D}}
\newcommand{\tA}{\mathbf{A}}
\newcommand{\Sym}[1]{\ensuremath{\mathcal{S}_{#1}}}
\newcommand{\simplexsymbol}[0]{\ensuremath{C}}

\newcommand{\R}[0]{\ensuremath{\mathbb{R}}}

\newcommand{\PFem}[1]{\ensuremath{\mathcal{P}_{#1}}}
\newcommand{\TheMethod}[0]{Geo\-de\-sic-BP}
\newcommand{\TheTitle}[0]{
Digital twinning of cardiac electrophysiology models from the surface ECG: A geodesic backpropagation approach}

\newcommand{\Vm}[0]{\ensuremath{V_m}}

\newcommand{\B}[0]{\ensuremath{B}} 
\newcommand{\V}[0]{\ensuremath{\mathcal{V}}}

\newcommand{\Gi}[0]{\ensuremath{\mathbf{G}_i}}
\newcommand{\Z}[1]{\ensuremath{Z_{#1}}}

\usepackage{cleveref}
\Crefname{equation}{Eq.}{Eqs.}

\usepackage[normalem]{ulem}
\usepackage{mdframed}
\newcommand{\revaat}[1]{#1}

\iftoggle{preprint}{
    
        \newcommand\copyrightnotice{%
        \begin{tikzpicture}[remember picture,overlay]
        \node[anchor=south,yshift=10pt] at (current page.south) {\fbox{\parbox{\dimexpr\textwidth-\fboxsep-\fboxrule\relax}{%
        \footnotesize 
\textcopyright 2023 IEEE.  Personal use of this material is permitted.  Permission from IEEE must be obtained for all other uses, in any current or future media, including reprinting/republishing this material for advertising or promotional purposes, creating new collective works, for resale or redistribution to servers or lists, or reuse of any copyrighted component of this work in other works. DOI: \href{https://dx.doi.org/10.1109/TBME.2023.3331876}{10.1109/TBME.2023.3331876}}}};
        \end{tikzpicture}%
    }
}{}

\iftoggle{preprint}{
    \title{\TheTitle{}}
    \author[1,4]{Thomas Grandits}
    \author[2]{Jan Verh\"ulsdonk}
    \author[1]{Gundolf Haase}
    \author[2]{Alexander Effland}
    \author[3]{Simone Pezzuto}
    \affil[1]{Department of Mathematics and Scientific Computing, University of Graz}
    \affil[2]{Institute for Applied Mathematics, University of Bonn}
    \affil[3]{Laboratory of Mathematics for Biology and Medicine, Department of Mathematics, University of Trento}
    \affil[4]{NAWI Graz, University of Graz}
}{
}

\date{}

\begin{document}

    \iftoggle{preprint}{
        \maketitle   
        \begin{abstract}
The eikonal equation has become an indispensable tool for modeling cardiac electrical activation accurately and efficiently. 
In principle, by matching clinically recorded and eikonal-based electrocardiograms (ECGs), it is possible to build patient-specific models of cardiac electrophysiology in a purely non-invasive manner. Nonetheless, the fitting procedure remains a challenging task.

The present study introduces a novel method, Geodesic-BP, to solve the inverse eikonal problem. 
Geodesic-BP is well-suited for GPU-accelerated machine learning frameworks, allowing us to optimize the parameters of the eikonal equation to reproduce a given ECG.

We show that Geodesic-BP can reconstruct a simulated cardiac activation with high accuracy in a synthetic test case, even in the presence of modeling inaccuracies. Furthermore, we apply our algorithm to a publicly available dataset of a \revaat{biventricular} rabbit model, with \revaat{promising} results.

Given the future shift towards personalized medicine, Geodesic-BP has the potential to help in future functionalizations of cardiac models meeting clinical time constraints while maintaining the physiological accuracy of state-of-the-art cardiac models.
        \end{abstract}    
        \copyrightnotice 
    }{
        \title{\TheTitle{}}
        \author{Thomas Grandits, Jan Verh\"ulsdonk, Gundolf Haase, Alexander Effland, and Simone Pezzuto
        \thanks{This work was supported by the Deutsche Forschungsgemeinschaft (DFG, German Research Foundation) through projects EXC-2047/1-390685813, and EXC2151-390873048.}
        \thanks{T.~Grandits and G.~Haase are with the Department of Mathematics and Scientific Computing, University of Graz, Austria (e-mail: \{thomas.grandits,gundolf.haase\}@uni-graz.at), T.~Grandits is additionally affiliated with NAWI Graz.}
        \thanks{J.~Verh\"ulsdonk and A.~Effland are with the In\-sti\-tu\-te for App\-lied Ma\-the\-ma\-tics, University of Bonn, Germany (e-mail: \{verhuelsdonk,effland\}@iam.uni-bonn.de).}
        \thanks{S.~Pezzuto is with the Laboratory of Mathematics for Biology and Medicine, Department of Mathematics, University of Trento, Italy,
        and with the Center for Computational Medicine in Cardiology, Euler Institute, Universit\`a della Svizzera italiana, Switzerland (e-mail: simone.pezzuto@unitn.it).}
        }

        \maketitle
        \begin{abstract}
            
        \end{abstract}  

    }

    \section{Introduction}

    \begin{figure*}[tb] 
        \centering
        \includegraphics[width=.9\linewidth]{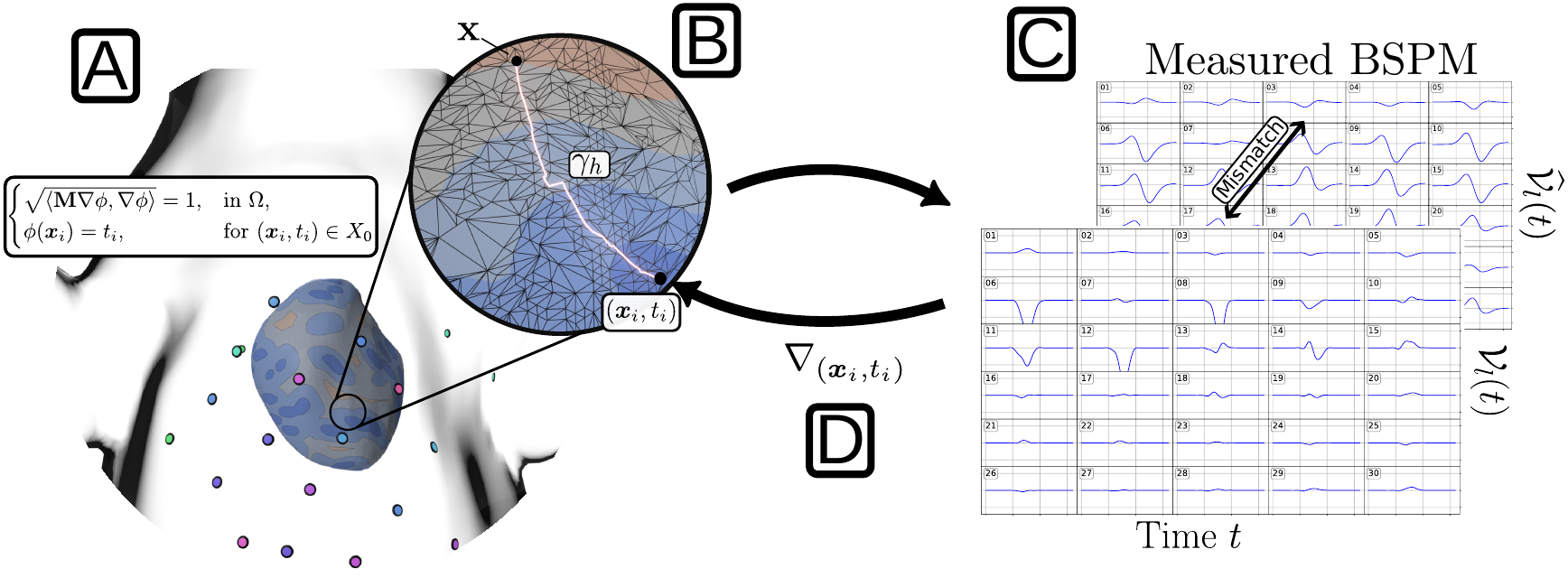}
        \caption{
            The anisotropic eikonal equation is an efficient tool to compute the electrical activation of a heart (A), here shown on a public rabbit torso model.
            Many solvers utilize the relation of the eikonal equation to geodesic paths, which in the case of \PFem{1}-meshes become piecewise-linear paths $\gamma_h$ between the initial conditions $\left( \vx_i, t_i \right)$ and points on the surfaces $\vx$ (B).
            By applying a non-linear transformation explained in \Cref{sec:ecg_method}, we can compute the measurable electrical activation on the torso surface and compute the mismatch against the measured potentials at sparse points (C).
            Implementing the approach in a differentiable manner using a backpropagation framework such as pytorch, we can efficiently calculate the gradient w.r.t.~the initial conditions and use it to optimally choose $\left( \vx_i, t_i \right)$ to match the measured potentials (D).
        }
        \label{fig:graphical_abstract}        
    \end{figure*}

    Cardiac digital twinning is a potential pillar of future's precision 
    cardiology~\cite{corral_acero_digital_2020}. A digital twin is a computational
    replica of the patient's heart, from its anatomy and structure up to its
    functional response. Commonly, it is obtained by fitting the parameters of
    a generic computational model to the clinical data of the patient. Here,
    we focus on cardiac electrophysiology, where we aim at personalizing a cardiac
    conduction model based on the eikonal equation, from non-invasive electrocardiographic (ECG) recordings~\cite{Pezzuto2021ECG,gillette_framework_2021}.

    Eikonal-based propagation models are an excellent basis for cardiac electrophysiology
    digital twins. They provide a good balance between computational efficiency and
    physiological accuracy for cardiac activation maps~\cite{keener_eikonal_curvature_1991,franzone_mathematical_2014,neic_efficient_2017}. When supplemented with the lead field approach, they can also provide accurate ECGs~\cite{pezzuto_evaluation_2017,gillette_framework_2021}. Hence, the parameters
    of eikonal-based models could be fitted from non-invasive clinical data---surface potentials and
    cardiac imaging---in an efficient manner.
    This goal can be achieved using optimization methods based on stochastic sampling strategies~\cite{camps_inference_2021,gillette_framework_2021}, Bayesian optimization~\cite{PezzutoBayes2022}, gradient-descent~\cite{grandits_geasi_2021}, or deep learning approaches~\cite{li_deep_2022}.

The numerical efficiency of the optimization procedure relies on two factors: the computational cost of the forward problem, and the number of samples required to achieve convergence within a given tolerance. Gradient-based methods can address the latter, if gradient computation is efficient and robust. For the former, many efficient eikonal solvers exist, but not all of them are suitable for automatic differentiation. A major criticality is that they internally solve a local optimization problem~\cite{jeong_fast_2008,grandits_fast_2021}, one for each element of the grid and for possibly many iterations. This may represent a bottleneck in the computation of the gradient. A possible solution is to compute the gradient for the continuous problem, which can be shown is equivalent to the computation of geodesic paths~\cite{grandits_geasi_2021, fehrenbach2023source}, and then discretize. The \emph{optimize-then-discretize} strategy is however not ideal, because it may introduce numerical noise in the gradient and degrade the convergence rate.

    In this paper, we introduce a novel inverse eikonal solver---named \emph{\TheMethod{}}---and use it to efficiently solve the ECG-based digital twinning problem \revaat{for ventricular models}.
    To this end, we start by introducing a highly-parallelizable anisotropic eikonal solver for triangulations in arbitrary dimension. The key contribution is a novel approach for the solution of the local optimization problem and based on a linear convex non-smooth optimization method.
    Locally, the linear geodesic path within each simplicial element is approximated with a guaranteed bound, similar to many presented eikonal solvers~\cite{jeong_fast_2008}.
    This solver can be implemented using current machine learning (ML) libraries to efficiently calculate the solution of the eikonal equation.
    The forward solver is suitable for the computation of the gradient via the backpropagation algorithm, effectively enabling a \emph{discretize-then-optimize} strategy (\TheMethod{}: Geodesic Back-Propagation).
    We will show that the backpropagation through such a solver is equivalent to piecewise-linear geodesic backtracking~\cite{grandits_inverse_2020,grandits_piemap_2021,grandits_geasi_2021}.
    Finally, we will demonstrate the potential of these methods for the problem of estimating initial eikonal conditions to match a given body surface potential map (BPSM) in a synthetic 2D study, and real 3D torso model.

In summary, the major contributions are as follows:
    \begin{enumerate}
        \item A novel, highly-parallelizable anisotropic eikonal solver, well-suited for GPGPU architectures and ML libraries.
        \item A backpropagation-based gradient computation, internally related to the computation of geodesics.
        \item An efficient implementation of the method.
        \item We demonstrate how to use \TheMethod{} for the identification of cardiac activation sequence from the non-invasive ECG recordings.
        \item We show that the inverse procedure is robust to model perturbations.
    \end{enumerate}
    
An outline of \TheMethod{} is provided in \Cref{fig:graphical_abstract}.

    \section{Material and methods}
    \label{sec:methods}
    In this section, the modeling assumptions for the electrical propagation inside the heart are introduced along with an efficient method for computing this activation in a parallelizable way.
    Then, we elaborate on the computation of the electrical activation given measurements on the torso leads, which ultimately leads to the formulation of the inverse ECG problem used throughout all experiments.

    \subsection{Cardiac activation}
    We model cardiac activation with the anisotropic eikonal equation (see~\cite{franzone_spreading_1993,pezzuto_evaluation_2017,neic_efficient_2017}), which reads as follows
    \begin{equation}
        \begin{cases}
            \sqrt{\left<\mathbf{M}(\vx) \nabla \phi(\vx), \nabla \phi(\vx)\right>} = 1, & \text{in $\Omega$},\\
            \phi(\vx_i) = t_i, & \text{for $(\vx_i, t_i) \in X_0$},
        \end{cases}
        \label{eq:anisotropic_eikonal_eq}
    \end{equation}
    where $\Omega \subset \R^d$, $d\ge 2$, is a bounded domain, $\phi: \Omega \to \R$ is the activation map, $\left<\cdot,\cdot\right>$ is the inner product in $\R^d$, and
    $\mathbf{M}: \Omega \to \Sym{d}$ is a symmetric positive definite (s.p.d.) tensor field, where $\Sym{d} \subset \R^{d\times d}$ denotes the space of $d$-times s.p.d. tensors, used to account for the orthotropic conduction velocity inside the heart~\cite{franzone_spreading_1993}.
    The onset of activation is determined by the set of initial conditions:
    $X_0 = \bigl\{ (\vx_i, t_i) \bigr\}_{i=1}^K$,
    where $\vx_i\in\Omega$ and $t_i$ are respectively the location and the onset timing of the $i$-th activation site. 
    For convenience, we additionally define the norm in the metric $\tD$, i.e. $\norm{\vx}_{\tD} = \sqrt{\left<\tD \vx, \vx\right>}$ for $\tD\in\Sym{d}$.

    Next, we discretize the domain $\Omega$ (assumed polytopal) by a simplicial triangulation. 
    We denote by $V = \{ \vv_i \}_{i=1}^{n_v}$ and $\mathcal{E} = \{ E_j \}_{j=1}^{n_e}$ respectively the set of vertices and elements of the triangulation. Each element
    $E_j\in\mathcal{E}$ is a $d$-simplex (triangle in 2-D, tetrahedron in 3-D) with $d+1$ vertices.
     We also define $E_{j\setminus i} \subset \partial E_j$ as the $(d-1)$-dimensional face of $E_j$ opposite to vertex $\vv_i$.
     Note that $E_{j\setminus i}$ is a $(d-1)$-simplex (triangle in 3-D, segment in 2-D).
    \revaat{For a general unstructured grid, we can} approximate the activation time with a piecewise linear finite element function $\phi: \Omega \to \R$, see~\cite{QV2008}. (For the sake of simplicity, we use the same symbol as for the continuous case.)
    For the nodal values we use the shorthand notation $\phi_i = \phi(\vv_i)$.
Following~\cite{jeong_fast_2008,fu_fast_2011,fu_fast_2013}, we approximate $\mathbf{D}(\vx) = \mathbf{M}^{-1}(\vx)$ as a piecewise-constant function and set $\mathbf{D}_j = \mathbf{D}(\vx)|_{E_j}$.

    We numerically solve the anisotropic eikonal equation \eqref{eq:anisotropic_eikonal_eq} by locally enforcing the following optimality condition for the nodal values:
    \begin{equation}
        \label{eq:hopflax_global}
            \phi_i = \min_{E_j \in \omega_i} \phi^*_{j\setminus i},
    \end{equation}
    where $\omega_i = \{ E \in \mathcal{E}: \text{$\vv_i$ is a vertex of $E$} \}$ is the set of elements containing $\vv_i$ and $\phi^*_{j\setminus i}$ is the minimum of the following optimization problem:
    \begin{equation}
    \label{eq:hopflax_local}
    \phi^*_{j\setminus i} = \min_{\vx \in E_{j \setminus i}} \phi(\vx) + \norm{\vv_i - \vx}_{\mathbf{D}_j}.
    \end{equation}
    The optimality condition~\eqref{eq:hopflax_global} is known as Hopf-Lax formula. It is possible to show that, under regularity assumptions on $\Omega$, $\mathbf{D}$, and \revaat{the triangulation}, the numerical solution of~\eqref{eq:hopflax_global} correctly convergences to the unique viscosity solution of~\eqref{eq:anisotropic_eikonal_eq} when the triangulation is refined~\cite{bornemann2006}.
    As explained later, this principle is at the basis of \TheMethod{}.

    \subsection{Solution of the local problem}
    \label{sec:methods_local}

    In what follows, we show how to efficiently and in a parallelizable way solve~\eqref{eq:hopflax_local}.
    We observe that~\cref{eq:hopflax_local,eq:hopflax_global} consists of a series of constrained minimization problems, one for each element of the patch. 
    Clearly, each local problem in~\eqref{eq:hopflax_local} is convex, because the objective function is convex in $\vx$ (note that $\phi$ is linear on $E_{j\setminus i}$) and the constraint is a simplex, thus convex as well.
    Instead of relying on exact formulas, e.g., available in 2-D and 3-D from~\cite{fu_fast_2013}, we propose here a new dimension-independent algorithm for the solution of~\eqref{eq:hopflax_local}.
    
\begin{figure}[tb]
        \centering
        \includegraphics[width=.23\textwidth]{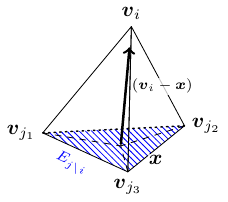}
        \includegraphics[width=.23\textwidth]{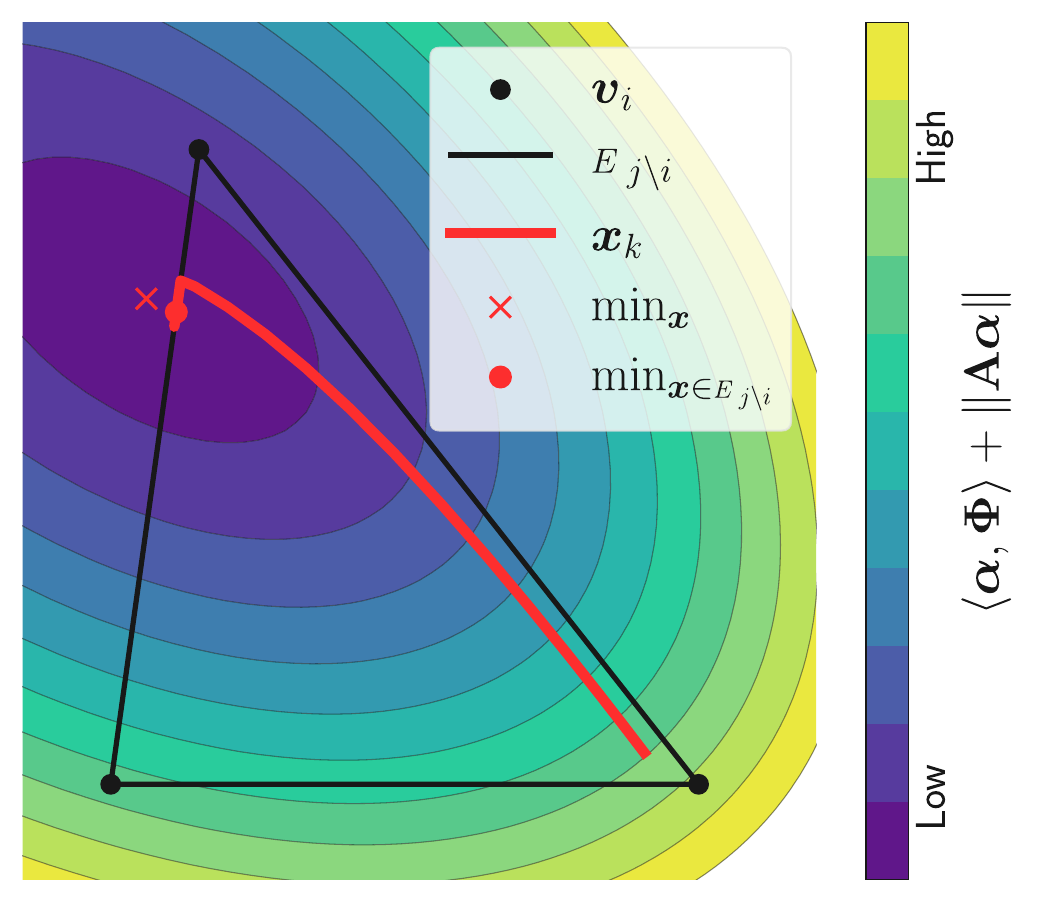}
        \caption{Conceptual visualization of a linear upwind discretization~\eqref{eq:upwind_eq_lsqr} on the example of a single tetrahedron (left).
        Finding the minimum upwind origin $\vx = \sum_n \alpha_n \vv_{j_n}$ can be cast into a non-smooth linear optimization problem, solvable using the FISTA algorithm (right).}
        \label{fig:upwind_update}
    \end{figure}

    With no loss of generality, we suppose that in~\eqref{eq:hopflax_local} the vertices of $E_j$ are $\{ \vv_{j_1}, \ldots, \vv_{j_{d+1}} \}$ and ordered such that the last vertex $\vv_{j_{d+1}}$ always corresponds to $\vv_i$. Next, we consider the set of barycentric coordinates for the face $E_{j\setminus i}$:
    \begin{equation}
        \simplexsymbol{}_{d} = \left\{ \balpha \in \R^{d}: \left< \balpha, \mathbf{1} \right> = 1 \land \balpha \geq 0 \right\},
        \label{eq:simplex_constraint}
    \end{equation}
    such that each point $\vx$ in a face $E_{j \setminus i}$ can be uniquely written as $\vx = \sum_{n=1}^{d} \alpha_n \vv_{j_n}$.
    Similarly, the \PFem{1} function $\phi(\vx)$ on $E_{j \setminus i}$ is given by the linear combination of the vertex values:
    $\phi(\vx) = \sum_{n=1}^{d} \alpha_n \phi_{j_n}$.
    Thus, \Cref{eq:hopflax_local} becomes
    \begin{equation}
    \begin{split}
        \phi^*_{j\setminus i} &= \min_{\balpha \in \simplexsymbol{}_{d}} \sum_{n=1}^{d} \alpha_n \phi_{j_n} + \Bigl\|\sum_{n=1}^{d} \vv_{j_n} \alpha_n - \vv_i \Bigr\|_{\tD_j} \\
        &= \min_{\balpha \in \simplexsymbol{}_{d}} \left< \balpha, \bPhi_{j \setminus i} \right> + \norm{\tA_{j \setminus i} \balpha}_2,
        \label{eq:upwind_eq_lsqr}
    \end{split}
    \end{equation}
    where $\bPhi_{j\setminus i} = \left( \phi_{j_1}, \ldots, \phi_{j_{d}} \right)^\top$ is the vector of nodal values of $\phi|_{E_j}$, except for $\phi_i$, and
    \begin{equation*}
        \mathbf{A}_{j \setminus i} = 
        \mathbf{D}_j^{1/2} \begin{pmatrix} 
            \vline &  & \vline \\
            \vv_{j_1} - \vv_{i} & \ldots & \vv_{j_d} - \vv_{i} \\
            \vline &  & \vline
        \end{pmatrix}.
    \end{equation*}
    A graphical representation of the local optimization problem is given in \Cref{fig:upwind_update}.

    \Cref{eq:upwind_eq_lsqr} now takes a common form of forward/backward-splitting algorithms, where the function to minimize is a linear combination of two convex functions with one potentially being non-smooth.
    Such problems can be minimized using the FISTA algorithm~\cite{beck_fast_2009}, which alternatingly optimizes using a proximal step into the non-smooth function and a gradient step in the direction of the smooth function.
    The smooth part of~\eqref{eq:upwind_eq_lsqr} was thus chosen as $h(\balpha) = \norm{\tA \balpha}_2$, and the non-smooth part as $g(\balpha) = \left< \balpha, \boldsymbol{\Phi} \right> + \delta_{\simplexsymbol{}_{d}}(\balpha)$, where $\delta_C$ is the indicator function on the set $C$ such that
    $$
    \delta_{C}(\vx) = 
    \begin{cases}
    0, & \text{if $\vx \in C$,}\\   
    \infty, & \text{otherwise.}
    \end{cases}
    $$
    The gradient of the smooth part $h(\balpha)$ is given by
    \begin{equation}
        \nabla h(\balpha) = \frac{\tA^\top\tA\balpha}{h(\balpha)},
        \label{eq:fista_norm}
    \end{equation}
    which is well-defined since $h(\balpha)>0$ for $\balpha\in C_d$ and non-collapsed elements.
    %
    The proximal step in the non-smooth part $g(\balpha)$ is then
    \begin{equation*}
        \begin{split}
        \operatorname{prox}_{\frac{1}{L}g}(\tilde{\balpha}) &= \argmin_{\balpha}  g(\balpha) + \frac{L}{2} \norm{\balpha - \tilde{\balpha}}^2 \\
        &=\proj_{\simplexsymbol{}_{d}} \left( \tilde{\balpha} - L^{-1} \bPhi \right),
        \end{split}
    \end{equation*}
    where $L$ is a constant.
    The projection $\proj_{\simplexsymbol{}_d}$ onto simplices is available in closed form~\cite[Figure 1]{duchi_efficient_2008}.
    The ISTA algorithm consists of the following update step:
    \begin{equation}
    \balpha_{k+1}=\proj_{\simplexsymbol{}_{d}} \Bigl( \left( \balpha_k - L^{-1} \nabla h(\balpha_k) \right)- L^{-1} \boldsymbol{\Phi} \Bigr).
    \end{equation}
    An optimal choice of $L$ in the algorithm is the Lipschitz constant of $\nabla h$, which we can bound as
    \begin{equation}
        L \le \frac{\norm{\tA^\top  \tA}}{\underline{h}}\left(1 + \frac{\norm{\tA}}{\underline{h}} \right),
        \label{eq:lipschitz_est}
    \end{equation}
    where $\norm{\cdot}$ in~\eqref{eq:lipschitz_est} refers to the spectral norm of the operator and $\underline{h} = \frac{1}{d \sqrt{d}} \sigma_{\min} \left(\tA\right) \le \norm{\tA \vx}_2 = h(\vx)$ is a lower bound of $h(\vx)$, for $\sigma_{\min} \left(\tA\right)$ being the smallest singular value of $\tA$.
    Note that $\tA$ is constant in each element, so the local Lipschitz constants can be efficiently pre-computed. 

    \begin{algorithm}[htb]
        \setstretch{1.0} 
        \caption{Solution of the local problem in \Cref{eq:upwind_eq_lsqr}}
        \label{alg:local_update}
        \KwData{$\tA$, $L$, $\bPhi = \left(\phi_1, \ldots, \phi_{d} \right)^\top$}
        \KwResult{$\balpha_{n_f}$} 
        $\hat{\balpha}_0 = \frac{1}{d} \mathbf{1}$ \\
        \For(\Comment*[h]{FISTA main loop}){$k = 1$ \KwTo $n_f$}{
            $h = \norm{\tA \hat{\balpha}_{k-1}}$ \\
            $\nabla h = h^{-1} \tA^\top \tA \hat{\balpha}_{k-1}$ \Comment{\Cref{eq:fista_norm}} 
            $\tilde{\balpha}_{k} = (\hat{\balpha}_{k-1} - L^{-1} \nabla h) - L^{-1} \bPhi$ \\
            $\balpha_{k} = \proj_{\simplexsymbol{}_{d}} (\tilde{\balpha}_k)$ \Comment{See \cite[Figure 1]{duchi_efficient_2008}}
            $\hat{\balpha}_{k} = \balpha_k + \beta_k (\balpha_k - \balpha_{k-1})$ \Comment{Acceleration} \label{ln:local_acceleration}
        }
    \end{algorithm}

    We summarize the solution of the local problem in \Cref{alg:local_update}, where we also include the acceleration on \Cref{ln:local_acceleration} with $\beta_k = \frac{k-1}{k+1}$, as per the FISTA method, see~\cite{beck_fast_2009}.
        Note that \Cref{alg:local_update} works for arbitrary spatial dimension $d\in \mathbb{N}$.     Since the local problem is convex, it has a unique minimum $\balpha^*$. 
    The rate of convergence is (see \cite[Theorem 4.4]{beck_fast_2009})
    \begin{equation}
        \phi^k_{j\setminus i} - \phi^*_{j\setminus i} \le \frac{2L \norm{\balpha_0 - \balpha^*}^2}{(k+1)^2}.
        \label{eq:fista_convergence}
    \end{equation}

    Finally, we compute the global eikonal solution by iteratively applying \Cref{alg:local_update} to each element $E_j \in E$ and for all vertices $\vv_i \in E_j$, and then take minimum: see \Cref{alg:fista_fim}.
    \iftoggle{singlecolalg}{
        \begin{algorithm}[htb]
            \setstretch{1.0} 
            \caption{\TheMethod{}}
            \label{alg:fista_fim}
            \KwData{$\mathcal{P}_1$ mesh $M = (V, \mathcal{E})$, 
            piecewise-constant metric $\tD$, initial conditions $X_0 = \bigl\{ (\vx_i, t_i) \bigr\}_{i=1}^K$}
            \KwResult{Geodesic distances $\phi_K$}
            $\forall \vv \in V: \phi_{-1}(\vv) = \phi_0(\vv) = \infty$ \\
            $\forall \left(\vx_i, t_i\right) \in X_0: \phi_{-1}(\vx_i) = \phi_0(\vx_i) = t_i$ \\
            \Comment{Compute all simplex upwind directions towards each vertex $i$ inside element $j$}
            $\tilde{E} = \bigcup_{i=1}^{d+1} \bigcup_j E_{j \setminus i}$ 
            Compute for each element $\tA_{j \setminus i}, L_{j \setminus i}$ \Comment{see \Cref{eq:upwind_eq_lsqr,eq:lipschitz_est}}
            \For{$k = 1$ \KwTo $K$}
            {
                $\phi_k := \phi_{k-1}$\\
                \Comment{Compute all updates in parallel}
                \Forpar($\tilde{E}_{j \setminus i} \in \tilde{E}$ \textbf{do}\label{ln:forpar}){ 
                    \Comment{\revaat{Collect} the simplex vertices opposite of the vertex $i$}
                    $(\vv_{j_1}, \ldots \vv_{j_d}) = \tilde{E}_{j \setminus i}$ \\
                    \Comment{Optional: Only compute updates when necessary}
                    \If{$\exists \vv_{j_n} \in \tilde{E}_{j \setminus i}: (\phi_{k-2}(\vv_{j_n}) - \phi_{k-1}(\vv_{j_n})) \ge \varepsilon$ \label{ln:necessary_mask}}
                    {
                        $\bPhi = (\phi_{k}(\vv_{j_1}), \ldots, \phi_{k}(\vv_{j_d}))^\top$\\
                        $\balpha_{n_f} =$ \Cref{alg:local_update}($\tA_{j \setminus i}, L_{j \setminus i}, \bPhi)$ \\
                        \With(\emph{enabled gradient computations})
                        {
                            \Comment{\Cref{eq:upwind_eq_lsqr}}
                            $\phi_{j\setminus i} = \left<\balpha_{n_f}, \bPhi \right> + \norm{\tA_{j \setminus i} \balpha_{n_f}} $ \label{ln:upwind_update} 
                        }
                    }
                }
                \With(\emph{enabled gradient computations})
                {
                    \Forpar($\vv_i \in V$ \textbf{do} ) {
                        \Comment{ see \Cref{eq:hopflax_global} \label{ln:global_psi_update}}
                        $\phi_k(\vv_i) = \min_{E_j \in \omega_i} \phi_{j\setminus i}$
                    }
                }
                \Comment{$\varepsilon$-convergence criterion}
                \If{$\forall \vv \in V: (\phi_{k-1}(\vv) - \phi_k(\vv)) < \varepsilon$ }
                {
                    $\phi_K = \phi_k$\\
                    \KwBreak 
                }
            }
        \end{algorithm}
    }{
        \begin{algorithm*}[htb]
            \setstretch{1.0} 
            \caption{\TheMethod{}}
            \label{alg:fista_fim}
            \KwData{$\mathcal{P}_1$ mesh $M = (V, \mathcal{E})$, piecewise-constant metric $\tD$, initial conditions $X_0 = \bigl\{ (\vx_i, t_i) \bigr\}_{i=1}^K$}
            \KwResult{Geodesic distances $\phi_K$}
            $\forall \vv \in V: \phi_{-1}(\vv) = \phi_0(\vv) = \infty$ \\
            $\forall \left(\vx_i, t_i\right) \in X_0: \phi_{-1}(\vx_i) = \phi_0(\vx_i) = t_i$ \\
            $\tilde{E} = \bigcup_{i=1}^{d+1} \bigcup_j E_{j \setminus i}$ \Comment{Compute all simplex upwind directions towards each vertex $i$ inside element $j$}
            Compute for each element $\tA_{j \setminus i}, L_{j \setminus i}$ \Comment{see \Cref{eq:upwind_eq_lsqr,eq:lipschitz_est}}
            \For{$k = 1$ \KwTo $K$}
            {
                $\phi_k := \phi_{k-1}$\\
                \Forpar($\tilde{E}_{j \setminus i} \in \tilde{E}$ \textbf{do} \Comment*[h]{Compute all updates in parallel \label{ln:forpar}}){ 
                    $(\vv_{j_1}, \ldots \vv_{j_d}) = \tilde{E}_{j \setminus i}$ \Comment{Contains the simplex vertices opposite of the vertex $i$}
                    \If(\Comment*[h]{Optional: Only compute updates when necessary} \label{ln:necessary_mask}){$\exists \vv_{j_n} \in \tilde{E}_{j \setminus i}: (\phi_{k-2}(\vv_{j_n}) - \phi_{k-1}(\vv_{j_n})) \ge \varepsilon$}
                    {
                        $\bPhi = (\phi_{k}(\vv_{j_1}), \ldots, \phi_{k}(\vv_{j_d}))^\top$\\
                        $\balpha_{n_f} =$ \Cref{alg:local_update}($\tA_{j \setminus i}, L_{j \setminus i}, \bPhi)$ \\
                        \With(\emph{enabled gradient computations})
                        {
                            $\phi_{j\setminus i} = \left<\balpha_{n_f}, \bPhi \right> + \norm{\tA_{j \setminus i} \balpha_{n_f}} $ \label{ln:upwind_update} \Comment{\Cref{eq:upwind_eq_lsqr}}
                        }
                    }
                }
                \With(\emph{enabled gradient computations})
                {
                    \Forpar($\vv_i \in V$ \textbf{do} ) {
                        $\phi_k(\vv_i) = \min_{E_j \in \omega_i} \phi_{j\setminus i}$
                    \Comment{ see \Cref{eq:hopflax_global} \label{ln:global_psi_update}}
                    }
                }
                \If{$\forall \vv \in V: (\phi_{k-1}(\vv) - \phi_k(\vv)) < \varepsilon$ }
                {
                    $\phi_K = \phi_k$\\
                    \KwBreak \Comment{$\varepsilon$-convergence criterion}
                }
            }
        \end{algorithm*}
    }
    
    \subsection{\TheMethod{}}
    \label{sec:backprop_method}
    %

	To utilize the computed eikonal solution in an inverse approach, we want to compute the gradients of the eikonal solution w.r.t.\ some quantity of interest, e.g., the initial conditions $X_0$.
    More precisely, we are interested in the gradient $\nabla_{\left(\vx_i, t_i\right)} \phi (\vx)$
    and its efficient computation.
    This can be done analytically for~\eqref{eq:anisotropic_eikonal_eq}, or algorithmically from \Cref{alg:fista_fim}. 
    The former (see e.g.~\cite{grandits_geasi_2021}) relies on a \emph{optimize-then-discretize} strategy. whereas the latter employs a \emph{discretize-then-optimize} approach. Following the latter, we compute $\nabla_{\left(\vx_i, t_i\right)} \phi (\vv)$ for all vertices $\vv \in V$ and observe that
    \begin{equation}
        \nabla_{\left(\vx_i, t_i\right)} \left(\phi (\vx) \right) 
        = \left< \balpha, \nabla_{\left(\vx_i, t_i\right)}  \left(\phi_{j_1}, \ldots, \phi_{j_{d+1}} \right)^\top \right>,
        \label{eq:dto_grad}
    \end{equation}
    for $\vx \in E_j$. Note that we neglect the derivative of $\balpha$ w.r.t.~$X_0$ as a consequence of the first variation of the geodesic distance \cite[Chapter 10]{oneill_semi_riemannian_1983}. 
    The nodal values $\phi_{j_n}$ in succession are then derived through~\eqref{eq:hopflax_local}, until we terminate at the simplex containing the initial condition $\vx_i \in E_k$.
    By defining $\vx_i$ in a continuous fashion (as was done in~\cite{grandits_geasi_2021}), such that for all vertices $\vv_{j_n}$ of the element $E_j$ containing $\vx_i$ it holds that
    \begin{equation}
        \phi_{j_n} = t_i + \norm{\vv_{j_n} - \vx_i}_{\mathbf{D}_j},
        \label{eq:upwind_xi}
    \end{equation}
    the derivative readily follows:
    \begin{equation*}
        \nabla_{\vx_i} \phi_{j_n} = -\frac{\tD_j \left( \vv_{j_n} - \vx_i\right)}{\norm{\vv_{j_n} - \vx_i}_{\tD_j}}, \quad \nabla_{t_i} \phi_{j_n} = 1,
    \end{equation*}
    %
    %
    which can be effectively calculated using automatic differentiation of~\eqref{eq:upwind_xi}, or subsequent nodes~\eqref{eq:dto_grad} which are connected by a geodesic path to $\left(\vx_i, t_i\right)$.
    In particular, we use the back-propagation algorithm implemented in many ML libraries. 
    More technical aspects are provided in \Cref{sec:practical_local}.

    \subsection{ECG computation}
    \label{sec:ecg_method}
    We describe the torso with a domain $\Omega_0 \subset \mathbb{R}^d$, whereas the heart is denoted by $\Omega$, that is assumed to be well contained in $\Omega_0$. The full heart-torso domain is $\Omega_T = \bar{\Omega}_0\cup\bar{\Omega}_T$, and the torso surface is $\Sigma = \partial\Omega_T$. The electrodes are placed on the torso at fixed locations $\tilde{X}_e = \left\{ \vx_{e_i} \right\}_{i=1}^N \subset \Sigma$. Each lead $\V_{l_e}(t)$ of the ECG is defined as the zero-sum of the electric potential at multiple locations. Thanks to the lead-field formula, we have the following representation:
    \begin{equation}
        \V_{l_e}(t)=\int_\Omega \left< \Gi(\vx) \nabla \Z{l_e}(\vx), \nabla \Vm(\vx,t) \right>\diff{\vx},
        \label{eq:integral_potential_ecg}
    \end{equation}
	where $\mathbf{G}_i : \Omega \to \Sym{d}$ is the intracellular conductivity tensor, $\Vm(\vx,t):\Omega\times [0,T] \to \mathbb{R}$ is the transmembrane potential, and $\Z{l_e} : \Omega_T \to \mathbb{R}$ is the $l_e$-th lead field, solution of the problem:
    \begin{equation}
    \left\{ \begin{aligned}
        & - \nabla \cdot \bigl( (\mathbf{G}_i + \mathbf{G}_e) \nabla \Z{l_e} \bigr) = 0, & \text{in $\Omega_T$,} \\
        & -\nabla\Z{l_e}\cdot\mathbf{n} = \delta_{\vx_{l_e}} - \frac{1}{|X_W|} \sum_{\vy\in X_W} \delta_{\vy}, & \text{on $\Sigma$.}
    \end{aligned} \right.
    \label{eq:lead_field}
    \end{equation}
    The set $X_W$ contains the electrodes at the left/right arm and left foot, used to form the Wilson's Central Terminal (WCT). \Cref{eq:lead_field} can be solved once, since the problem is time-independent. The bulk conductivity of the torso, $\mathbf{G}_i + \mathbf{G}_e$, also depends on the extracellular conductivity tensor $\mathbf{G}_e: \Omega \to \Sym{d}$. The derivation of this problem and exact formulation of~\eqref{eq:lead_field} can be found in further detail in~\cite{franzone_mathematical_2014,keener_mathematical_1998,grandits_geasi_2021}.

    
    Next, we assume that the transmembrane potential is travelling wave with activation times dictated by $\phi(\vx)$:
    \begin{equation}
        \Vm(\vx,t) = U(t-\phi(\vx)) = U(\xi).
    \end{equation}
     Since we are only interested in the depolarization sequence, we use the following template action potential, see~\cite{grandits_geasi_2021,pezzuto_evaluation_2017}:
    \begin{equation}
        U(\xi) = K_0 + \frac{K_1 - K_0}{2} \tanh \left( 2 \xi/\tau \right),
        \label{eq:tanh_waveform}
    \end{equation}
    where $K_0$ and $K_1$ are respectively the minimum and maximum transmembrane potential.

    Additionally, we precompute the discretized integral ECG operator~\eqref{eq:integral_potential_ecg} by means of a 2nd order simplicial Gaussian quadrature scheme.
    We used the library scikit-fem~\cite{skfem_2020} to assemble this operator for our \PFem{1} mesh, automatically interpolating and integrating $\Vm$ at the Gaussian quadrature points.
    This way, we can compute the leads as $\V_{l_e}(t) = \left(\B U(t - \phi)\right)_{l_e}$ where $\B: \R^{n_v \times N}$ is the mentioned discretized, precomputed lead field integration operator of~\eqref{eq:integral_potential_ecg}, $n_v = \abs{V}$ are the number of vertices and $N$ is the number of leads.
    
    For optimizing the inverse ECG problem, we want to fit $\V_{l_e}$ of~\eqref{eq:integral_potential_ecg} against a measured $\hat{\V}_{l_e}$.
    We consider for that purpose the minimization problem
    \begin{equation}
        \begin{split}
        \min_{(\vx_i, t_i) \in X_0} &\frac{1}{N \abs{T}} \sum_{l_e} \int_T \left(\V_{l_e}(\phi_{X_0}, t) - \hat{\V}_{l_e}(t) \right)^2 \diff{t},
        \end{split}
        \label{eq:inverse_ecg}
    \end{equation}
    where $\phi_{X_0}$ is the eikonal solution for given initial conditions and $N$ defines the number of measured/computed leads.
    Note that with $\B$ computed from~\eqref{eq:integral_potential_ecg}, we can easily conclude for  $\xi = t - \phi$ that $\nabla_{\left(\vx_i, t_i\right)} \V_{l_e}(t) = -\B^\top \frac{\partial U(\xi)}{\partial \xi} \nabla_{\left( \vx_i, t_i \right)} \phi$, where $\nabla_{\left( \vx_i, t_i \right)} \phi$ can be computed using \Cref{alg:fista_fim} as discussed in \Cref{sec:backprop_method}.
    Note that backpropagation through $\nabla_{\left(\vx_i, t_i\right)} \V_{l_e}$ will yield the equivalent result.
    The resulting gradients were then used for a gradient-based optimization utilizing ADAM~\cite{kingma_adam_2017}.


	\subsection{Implementation aspects}
    \label{sec:practical_local}
        There are a few practical remarks when implementing \Cref{alg:local_update,alg:fista_fim} that are important to consider:
            When computing the local update in \Cref{alg:local_update}, the number of iterations $n_f$ is the dominating factor for performance. Since the local problem is strictly convex with a precomputed Lipschitz constant $L_E$, the exact error bounds (see \eqref{eq:fista_convergence}) towards the global minimum in~\eqref{eq:upwind_eq_lsqr} are known and could be used to compute $n_f$ on an element basis. 
            Note that $L_j$ only depends on $\mathbf{A}_{j\setminus i}$, which is strongly related 
            to the mesh geometry and may increase as element quality decreases.
            
            Furthermore, The upper bound presented in~\eqref{eq:lipschitz_est} is not tight and might lead to slow convergence of the local FISTA updates in \Cref{alg:local_update}.
            In practice, we thus recomputed $\underline{h}$ by minimizing $\balpha^*$ through \Cref{alg:local_update} with $\bPhi = \mathbf{0}$ and using the estimate $\underline{h} = h(\balpha^*) / 1.1$.


            When considering \Cref{alg:fista_fim}, note that the parallelization denoted in \cref{ln:forpar} is more descriptive for CPU-parallel implementations. In ML frameworks, in contrast, the parallelization comes from the vectorization of single implemented operations, which might require re-ordering of some operations. 
            Still, all operations inside the loop starting from \cref{ln:forpar} can be parallelized, where only a special parallelization is required for \Cref{ln:global_psi_update} (scatter operation).
            This makes it well-suited for vectorized implementations on a GPU in frameworks such as PyTorch~\cite{paszke_pytorch_2019}.

            The optional performance optimization denoted in \cref{ln:necessary_mask} already significantly cuts down on the required memory for backpropagation, only computing updates when the solutions of simplex elements have changed.
            Checking for an $\varepsilon$ change in activations before the update can also avoid unwanted gradient eliminations in cases where the old and new values of $\phi$ take the same value.

            Note that the local problem in \Cref{alg:local_update} could also be solved using the exact formulation derived for low-dimensional cases $(d \le 4)$ and fixed simplex dimensions as in previous works~\cite{jeong_fast_2008,fu_fast_2013}.
            In practice, however, the parallelization of this explicit solution is limited by the required condition to satisfy the constraint in~\eqref{eq:simplex_constraint}.
            Such branching leads to unwanted performance losses in SIMD architectures when parallelizing over multiple elements.

            \revaat{ 
                Note that the number of activation sites is never modified by \Cref{alg:fista_fim}. Some activation sites may however be overshadowed by the surrounding activation, thus effectively providing no contribution to the overall activation. 
                When this happens, the gradient w.r.t.\ the overshadowed site is simply zero, thus the point is ignored by the algorithm.
                If an initial condition $(\vx, t_i)$ is overshadowed as described at the final iteration of the optimization, we often refer to it as inactive.
            }



    \section{Results}
    \label{sec:results}
    In this section, details of the dataset and results of applying our method for matching an ECG for both a synthetic 2-D and a real 3-D case are provided.



    \subsection{Simulation study}
    \label{sec:idealized_torso}
    \subsubsection{Setup}
\begin{figure}[htb]
        \centering
        \includegraphics[width=.7\linewidth]{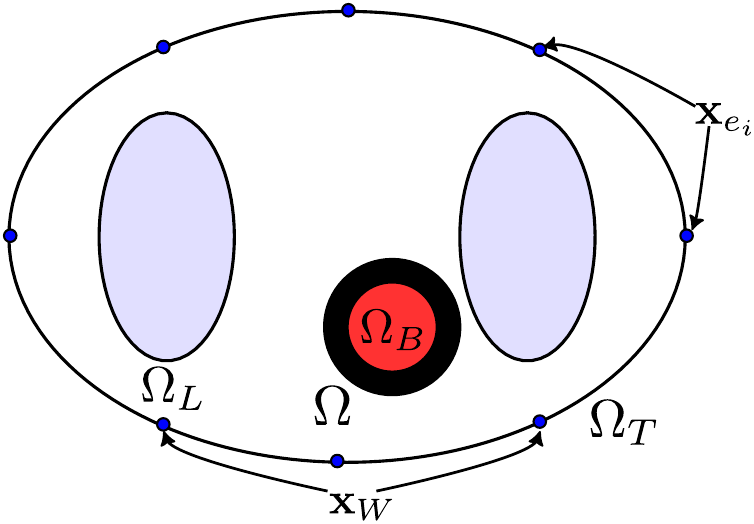}
        \caption{The 2D simulation setup.
        We consider a simple symmetrical left ventricular slice $\Omega$ to be located in a heterogeneous torso $\Omega_T$, containing the blood pool $\Omega_B$ and lungs $\Omega_L$, with different conductivities.
        $8$ electrodes $\vx_{e_i}$ (blue circles) are circularly distributed around the torso, two of which are chosen as the WCT $\vx_W \in X_W$.}
        \label{fig:idealized_torso_setup}
    \end{figure}

    To test and analyze our algorithm, we first start with a simple 2-D simulation study. 
    We designed a simple setup consisting of an idealized torso with $8$ electrodes distributed around the surface, with left and right arm electrodes composing Wilson's central terminal (WCT) (see \Cref{fig:idealized_torso_setup}). 
    Seven lead fields were then computed each w.r.t.~WCT.

    \begin{figure*}[tbh]
        \centering
        \includegraphics[width=.9\linewidth]{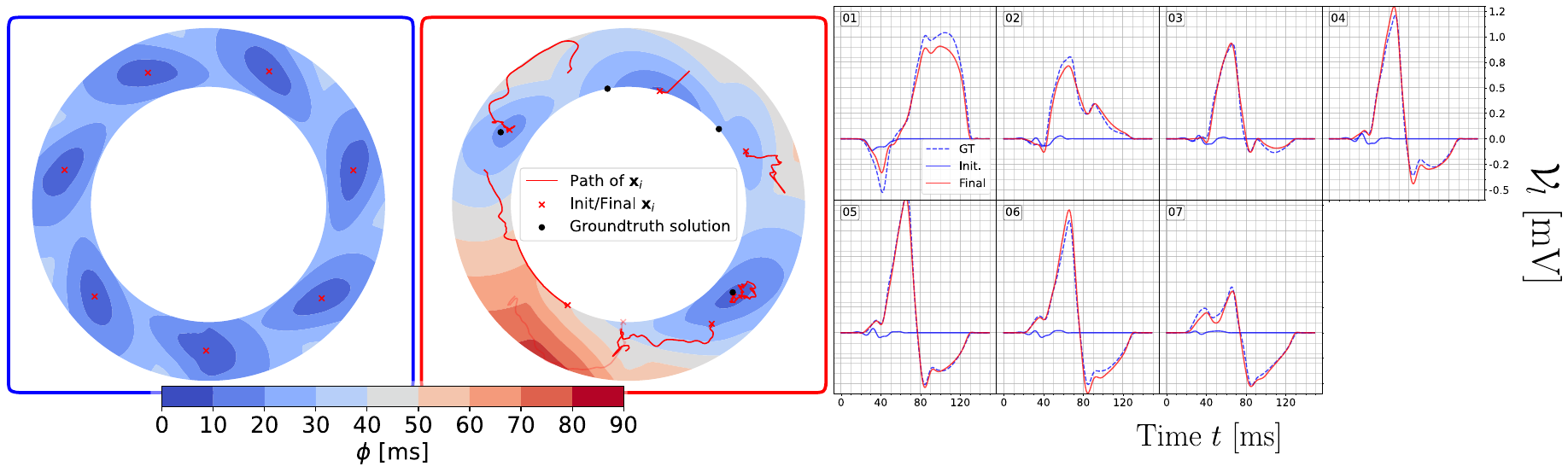}
        \caption{Result of the ECG optimization on the idealized torso experiment using $7$ leads.
        The initial model (left) uses a simple simultaneous activation, evenly spread throughout the myocardium.
        The optimization reduces the number of initial conditions and is able to match both $\phi$ and $\V_l$ closely (middle and right).
        The transparent red lines in the middle plot show the optimization path of initial conditions $(\vx_i, t_i)$, which are \revaat{inactive} in the final solution.}
        \label{fig:idealized_torso_results}
    \end{figure*}
    
    The LV, meshed at a resolution of \SI{0.9}{\milli \meter} ($\approx$ \num{2700} DoFs), contained 4 initial condition tuples $\left(\vx_i, t_i \right)$, located throughout half of the LV myocardium with timings between \SI{2.5}{\milli \second} and \SI{10}{\milli \second}.
    For computing the ground truth, we use another model with a \SI{0.45}{\milli \meter} space resolution and perturbed conductivities in all sub-domains of the torso.
    As an initial guess, we consider 8 activation sites, all with \SI{0}{\ms} timing and evenly distributed at the transmural center.
    The optimization was performed using ADAM~\cite{kingma_adam_2017} for 400 epochs using a learning rate of \num{0.5}, taking $\approx$ \SI{5}{\minute} on a AMD Ryzen 7 3800X 8-core processor.

    \subsubsection{Results}
    The results of the optimization against the initial and ground-truth solution and parameters are shown in \Cref{fig:idealized_torso_results}.
    We can nearly perfectly fit the ECG and approximately recover the initial conditions.
    Overall, the total root mean squared error (RMSE) w.r.t.~the eikonal solution $\phi$ is only \SI{5.8}{\milli\second}.
    We tested the algorithm also for a varying number of initial conditions and found that it is a hyperparameter of importance, but need not be chosen exactly:
    A cross-validation on the 2-D torso revealed that the worst match is achieved, by underestimating the number of initial conditions.
    In such a case, the achievable activations are not complex enough and the model is unable to faithfully replicate both the ECG and unobservable eikonal solution.
    Increasing the number of initial conditions beyond the ground truth's number is neither able to better match the ECG by a significant margin, but will also not overfit and significantly increase the error on the eikonal solution.
    The most common cause of high errors on the eikonal solution is related to getting trapped in local minima when optimizing for the ECG, circumventable by batched optimization (not utilized in any of the shown experiments).
    
    \subsection{Body Potential Surface Map Reconstruction}
    \label{sec:bpsm_exp}
    \subsubsection{Setup}
    \label{sec:data_used}
    \color{black} 
    The second anatomy is from a public rabbit ECGi model~\cite{moss_volumetric_2022}.
    The anesthetized rabbit was equipped with a 32-lead ECG vest during normal sinus rhythm and the signals were recorded with a frequency of \SI{2048}{\hertz}.
    The ECG was recorded for 10 seconds, filtered, and averaged over all beats to produce a single beat.
    From this beat, we extracted the QRS duration to match the biventricular activation.
    The biventricular mesh consisted of \num{508 975} DoFs and \num{2 522 501} tetrahedra, roughly equivalent to a mean edge length of \SI{0.82}{\milli \meter} in a human heart (assuming a size factor of 3 between rabbit and human hearts).
    Note that we also tested the optimization with a lower resolution of \num{143 564} DoFs and \num{591 975} tetrahedra, roughly equivalent to a mean edge length of \SI{1.25}{\milli \meter}.

    We use the electrophysiological parameters from~\cite{moss_computational_2022}, including cardiac fiber orientation, lead fields $\Z{l}$, and conduction velocity.
    In the conducted experiment, we chose $K_0 = \SI{-85}{\milli\volt}$, $K_1 = \SI{30}{\milli\volt}$, $\tau_1 = \SI{1}{ms}$ in~\eqref{eq:tanh_waveform}.
    The extracted QRS sequence was chosen during the time interval $t \in \SIrange{100}{165}{\ms}$.
    We optimized for the initial conditions $(\vx_i, t_i)$ of the model, i.e.~$X_0$ in~\Cref{alg:fista_fim}.
    For this purpose, we randomly distributed 100 points on the biventricular surface ($\vx_i \in \partial \Omega $) and with all timings initialized to the approximate Q onset $t_i = \SI{20}{\milli \second}$.
    
    \subsubsection{Results}

    \label{sec:runtime}
    The experiment was run for 400 epochs using a learning rate of \num{0.5} on an Nvidia A100 GPU (40GB), taking approximately \SI{1.8}{\hour} (or \SI{0.7}{\hour} on the reduced resolution).
    The computed results can be found on Zenodo~\cite{grandits_geodesic_bp_dataset_2023}.

    \begin{figure*}[htb]
        \centering
        \includegraphics[width=.9\linewidth, trim={0cm 0.25cm 0cm 0.25cm}, clip]{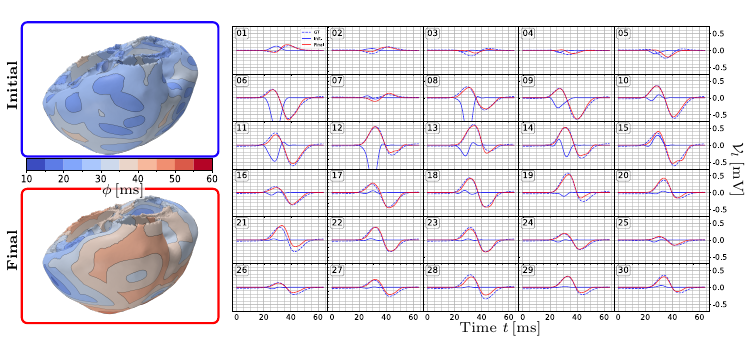}
        \caption{Result of the ECG optimization on the rabbit heart model.
        We compare the initial and final model (blue/red) and show the computed activations $\phi(\vx)$ (left) and BPSM $\V_l(t)$ for each model (right).
        The computed ECGs (solid lines) are shown on top of the measured ECGs (dashed blue line). 
        We show the first 30 of 31 available leads to save an extra row.}
        \label{fig:bpsm_result1}
    \end{figure*}
    The results for optimizing~\eqref{eq:inverse_ecg} with 400 epochs of ADAM using \Cref{alg:fista_fim} can be seen in \Cref{fig:bpsm_result1}.
    The intial and optimized activation are shown in the left upper and lower side of the image. is shown together with the comparison of all ECGs (right).
    The initially randomly sampled points evenly cover the domain $\Omega$, to allow for a fast excitation of the whole domain, though this leads to overexpressed deflections in many leads initially.
    The optimization significantly changes the activation of the heart and is able to achieve a promising match of the ECG in many leads.
    The loss in~\eqref{eq:inverse_ecg} converges after $\approx$ 250 epochs, with only minor variations in the parameters thereafter.
    Of the initial 100 points, only 91 are still active at the end of the optimization, meaning that the corresponding $t_i$ is contained in the final solution and not earlier activated by closer initial conditions.

    \section{Discussion}
    We presented \TheMethod{}, an algorithm to efficiently solve the inverse ECG problem while being very efficient on a single GPU. We are able to fit the parameters of a patient-specific model only from the surface ECG and the segmented anatomy. The fitted model faithfully reproduced the recorded ECG in an animal experiment. 
    In a synthetic case \revaat{(\Cref{sec:idealized_torso})}, \TheMethod{} is able to recover the ground-truth solution with high accuracy, even in the presence of perturbation in the forward model and in the lead field.
    The resolution of the considered 3-D heart model consists of $>10^5$ DoFs, which is sufficiently accurate for a high-fidelity eikonal solution in a human heart, as it would correspond to roughly \SI{0.8}{\milli \meter} edge length. 
    At full resolution, all computations only take \SI{16}{\second}.
    

    The considered cardiac parameter estimation problem is of high-relevance~\cite{ruiz_herrera_physics_informed_2022,coveney_probabilistic_2020,corral_acero_digital_2020,roney_technique_2019}, yet only a few works are able to efficiently solve this problem using solely non-invasive ECG measurements.
    The considered heart model and ECG originates from a real-world rabbit rather than a human model since it is (to the knowledge of the authors) the only publicly available full torso and heart model, captured with this high-level detail and including preprocessed electrophysiological measurements.   
    Such animal models share many properties of human cardiac electrophysiology~\cite{odening_animals_2020}.
    Furthermore, the algorithm and equations involved translate 1:1 to human cardiac electrophysiology.
    \revaat{Likewise, the methods involved are not restricted to biventricular models, and can be extended to atrial, or whole heart models under sinus rhythm.}
    \revaat{We also foresee the future possibility of functionalizing publicly available human heart and torso models, currently missing parameterization~\cite{aras_experimental_2015}.}

    The runtime in the full rabbit heart model of \SI{1.8}{\hour} for the computed 400 iterations is still expensive but it reduces to \SI{0.7}{\hour} for a coarser resolution (see \Cref{sec:runtime}). 
    \revaat{
    Two main advantages of our novel method are the high number of parameters that can be simultaneously optimized, allowing for many initial sites in contrast to other works~\cite{gillette_framework_2021,Pezzuto2021ECG}, and having no spatial restrictions on the onset location~\cite{PezzutoBayes2022}.
    The first advantage is a consequence of utilizing backpropagation to compute the gradient w.r.t.~our parameters, resulting in the runtime not being dependent on the number of parameters 
    and the convergence being dictated by properties of the underlying function to optimize~\cite{beck_first_order_2017}.
    Thus, compared to existing works based on derivative-free optimization with the same objective, our algorithm will scale much better in terms of spatially-varying parameters such as conduction velocities~\cite{conn_introduction_2009}.} 
    \revaat{
    The latter advantage of arbitrary onset locations follows from the use of the volumetric description of the eikonal equation and lead fields, from which we can compute parameters and the local activation times not only on the heart surface, but also transmurally.
    In conjunction, both of these advantages allow us to model potentially more complex activations than previously were able by other works.

    From an implementation point of view, \Cref{alg:fista_fim} can be readily applied to optimizing the conduction velocity tensor $\mathbf{D}$ with little modification.
    Based on our previous research~\cite{ruiz_herrera_physics_informed_2022,grandits_piemap_2021}, however, optimizing conduction velocities either requires strong prior information in the form of regularization, or multiple activation sequences to create meaningful results.
    The question how well such information can be reconstructed from the BSPM, and the space of possible solutions w.r.t.~initial conditions both remain topics for future research.

    One major limitation of \TheMethod{} is that it cannot handle re-entry phenomena such as tachycardia and fibrillation, because it relies on the eikonal equation with focal boundary conditions. However, the eikonal model can be reformulated to allow re-entry, see~\cite{pernod2011multi} and it could be used in an optimization problem where the initial condition is to be identified from the ECG.
    }
    %
    Finally, we also foresee many practical improvements for reducing the computational time of the local solver as outlined
    \Cref{sec:practical_local}.

\section{Conclusions}
    This paper presents \TheMethod{}, an efficient method to compute piecewise linear geodesics, and demonstrated its performance for solving the inverse ECG problem in a limited parameter setting.
    By creating the link between backpropagation through a FEM-solver with geodesic backtracking, we were able to show that our previous FEM-based inverse eikonal models~\cite{grandits_inverse_2020,grandits_piemap_2021} are linked to our recent work~\cite{grandits_geasi_2021} working with geodesic backtracking.
    Such inverse procedures are current research topics of high-relevance and efficient inverse eikonal models such as the present one may help in advancing patient-specific cardiac electrophysiology to practical real-world scenarios in the near future.

    \iftoggle{preprint}{
        \printbibliography 
    }{
        \bibliographystyle{ieeetr}
        \bibliography{literature}
    }

\end{document}